\documentclass[10pt,a4paper,draft]{article}
\usepackage{amsmath}
\usepackage{amsfonts}
\usepackage{amssymb}
\usepackage{amscd}
\makeatletter\@addtoreset{equation}{section}\makeatother
\newtheorem{theorem}{Theorem}[section]

\newtheorem{lemma}{Lemma}
\def\Z{{\mathbb Z}}
\def\F{{\mathbb F}}

\newtheorem{definition}{Definition}

\begin{document}
\title{{\bf
Indecomposable set-theoretical solutions to the Quantum Yang-Baxter
Equation on a set with prime number of elements}}
\author{ {\Large Pavel Etingof, Robert Guralnick, Alexander Soloviev}\\
\date{}
}
\maketitle
\section{Introduction}
In this paper we show that all indecomposable nondegenerate 
set-theoretical solutions to the Quantum
Yang-Baxter equation on a set of prime order 
are affine, which allows us to give 
a complete and very simple classification of such solutions. 
This result is a natural application of
the general 
theory of set-theoretical solutions to the quantum Yang-Baxter equation,
developed in \cite{ESS},\cite{LYZ},\cite{S} following a suggestion of 
Drinfeld \cite{Dr}. It is also a 
generalization of
the corresponding statement for involutive set-theoretical solutions proved
in \cite{ESS}. 

In order to prove our main result, we use the theory 
developed in \cite{S} to reduce the problem to a group-theoretical 
statement: a finite group with trivial center 
generated by a conjugacy class of prime order is a subgroup of the affine 
group. Unfortunalely, we did not find an elementary proof of 
this statement, and our proof relies on the classification
of outer automorphisms of finite simple groups. 

The structure of the paper is as follows. In Section 2 we give background
material and formulate key theorems about set-theoretical solutions
to the QYBE.  In Section 3 we prove
the main theorem. In the appendix 
we prove the above group-theoretical 
statement, used in the proof of the main theorem. 

\section{Set-theoretical solutions to the Quantum Yang-Baxter Equation}
\subsection{Structure group and group-theoretical characterization of
nondegenerate braided sets}
Let $X$ be a nonempty set and $S:X\times X\to X\times X$ a
bijective map. We call a pair $(X,\ S)$ a braided set if the following
braiding condition holds in $X\times X\times X$:

\begin{equation}
\label{braid}
S_1S_2S_1=S_2S_1S_2,
\end{equation}
where $S_1=S\times id$, $S_2=id\times S$.

{\noindent \underline{Remark.}} Consider the map
$R:X\times X\to X\times X$ given by $R=\sigma S$, where
$\sigma(x,y)=(y,x)$ for $x,\ y\in X$. Then $(X, S)$ is a braided set
if and only if $R$ satisfies the Quantum Yang-Baxter equation. 

We also introduce the
maps $g:X\times X\to X$ and $f:X\times X\to X$ as
components of $S$, i.e. for $x,\ y\in X$
 $$S(x,y)=(g_x(y),f_y(x)).$$
\begin{definition}
(i) We call a set $(X,S)$ nondegenerate if $g_x(y)$
is a bijective function of $y$ for fixed $x$ and $f_y(x)$ is a bijective
function of
$x$ for fixed $y$.
(ii) We call a  set $(X,S)$ involutive if $S^2=id_{X^2}$.
\end{definition}

In \cite{ESS},\cite{LYZ},\cite{S} authors developed
a theory of nondegenerate braided sets and gave a
description (cf. \cite{S}) of
the category of such sets in group-theoretical terms. We will be using
 this theory in our paper and will formulate 
the necessary results along the
way.

 From now we always assume $(X, S)$ to be a nondegenerate braided set and
refer to it as a "solution" keeping in mind that $\sigma S$ is a
solution to the QYBE.

It is useful to associate with a solution $(X,S)$ two
groups $G_X$ and $A_X$.
\begin{definition}
Define the group $G_X$ as the group generated
by the elements of $X$ subject to the relations $xy=y_1x_1$
if $S(x,y)=(y_1,x_1)$, where  $x, y \in
X$. We call $G_X$ the structure group of the solution $(X,\ S)$.
\end{definition}
\begin{definition}
Define the group $A_X$
as the group  generated by the elements of $X$ subject to relations
$x_1\bullet y=y_2\bullet x_1$,  where $x, y \in X$ and $x_1,\ y_2\in
X$ are defined by $S(x,y)=(y_1,x_1),\
S(y_1,x_1)=(x_2,y_2)$. We
call $A_X$ the derived structure group of the solution $(X,\ S)$.
\end{definition}

Since $G_X$ and $A_X$ are generated by $X$, 
there are natural maps $\psi_G:X\to G_X$
and $\psi_A:X\to A_X$. 

\begin{definition}
A solution $(X,S)$  is called injective if
the map $\psi_G:X \to G_X$ is injective.
\end{definition}

\begin{theorem}\cite{S}\label{invo}
A solution $(X, S)$ is involutive if and only if it is
injective and
its derived structure  group $A_X$ is abelian.
\end{theorem}

The following theorem is the first step towards establishing a bridge
between solutions and their counterparts in the group
world: bijective cocycle 7-tuples.

\begin{theorem}
\label{action}
(i) The map $x\to f_x^{-1}$ can be extended to a  left action of $G_X$ on
$A_X$ by automorphisms. 

(ii) The map $\bar{\rho}:X\times X\to Aut(X)$ given by the formula
$\bar{\rho}(x,y)(z)=f_x^{-1}(f_y(g_{f_z^{-1}(y)}(z)))$ can be extended
to the action $\rho:G_X\ltimes A_X\to Aut(X)$ of 
the semidirect product $G_X\ltimes A_X$ on
$X$ such that $\bar{\rho}=\rho(\psi_G\times \psi_A)$.
\end{theorem} 
{\bf \noindent Proof:}

The statements of the theorem easily follow from Theorems 2.4,2.7
in \cite{S}.
$\blacksquare$

Since the goal of this paper is to study solutions of prime order, 
we will assume from now on that the set $X$ is {\it finite}. 

Let $(X, S)$ be a solution.
\begin{theorem}
\label{Gamma}
There are a $G_X$ - invariant central subgroup $\Gamma_2\subset A_X$
of finite index in $A_X$ and a normal subgroup  $\Gamma_1\subset G_X$ such
that
(i) $\Gamma_1$ acts trivially on $X$ and thus on $A_X$ yielding an action
of $G_X/\Gamma_1$ on $A_X/\Gamma_2$, $\rho_{\Gamma}:G_X/\Gamma_1\to
Aut(A_X/\Gamma_2)$.

(ii) There is a bijective 1-cocycle 
$\bar{\pi}:G_X/\Gamma_1\to A_X/\Gamma_2$ with respect to action
$\rho_{\Gamma}$, i.e. $\bar{\pi}$ satisfies the relation
$\bar{\pi}(ab)=\rho_{\Gamma}(b^{-1})(\bar{\pi}(a))\bar{\pi}(b)$.
\end{theorem}

The action $\rho:G_X\ltimes A_X\to Aut(X)$ was instrumental (cf. \cite{S})
for
establishing a 1-1 correspondence between solutions and bijective cocycle
7-tuples. It is crucial here as well since it allows us to understand
the indecomposability property.

\subsection{Indecomposable solutions}
\begin{definition}
We call a solution $(X, S)$ decomposable if there is a
partition of $X$ into two disjoint nonempty subsets $X_1,
X_2$ such
that 
$S(X_1\times X_1)\subset X_1\times X_1$, and  
$S(X_2\times X_2)\subset X_2\times X_2$. 
\end{definition}

It is clear that in this case $(X_1, S|_{X_1\times X_1})$, 
$(X_2, S|_{X_2\times X_2})$ are also solutions.

\begin{definition}
We call a solution $(X, S)$ indecomposable if it is not decomposable.
\end{definition}

The following lemma plays a key 
role in studying indecomposable solutions.
\begin{lemma}
A solution $(X, S)$ is indecomposable if and only if
the action $\rho$ of Theorem \ref{action} is transitive. 
\end{lemma} 
{\bf \noindent Proof:} If $(X, S)$ is decomposable into 
$(X_1, S|_{X_1\times X_1})$ and $(X_2, S|_{X_2\times X_2})$ then due to
nondegeneracy and finiteness of $(X, S)$ 
one has $S(X_1\times X_2)=X_2\times X_1$
and $S(X_2\times X_1)=X_1\times X_2$. This implies that $X_1$, $X_2$ are
orbits for $\rho$. Similarly, one can verify that if $X_1$, $X_2$ are
nonempty $\rho$-orbits such that $X=X_1\cup X_2$ then $(X, S)$ is
decomposed into nondegenerate subsolutions $(X_1, S|_{X_1\times X_1})$,
$(X_2, S|_{X_2\times X_2})$.$\blacksquare$ 

\begin{definition}
A solution $(X, S)$ of the form $S(x,y)=(\phi(y,x),x)$ for some
$\phi:X\times X\to X$ is called derived.
\end{definition}

{\bf\noindent Example (see \cite{Dr},\cite{LYZ},\cite{S}}\label{der})
Let $G$ be a group acting on itself by conjugation and by $\rho_c$ on a
set $X$ in such a way that a map $i:X\to G$ is equivariant. Then $(X, S)$
with $S(x,y)=(\rho_c(i(x))(y),x)$ is a derived solution.

Results in \cite{S} imply the following:

\begin{theorem}
\label{derived} Let $(X,S)$ be a derived solution. 

(i) $(X,S)$ is isomorphic to  the solution described in Example \ref{der}
for the group $G=G_X/\Gamma_1$ and the map $i$ coming from $\psi_G:X\to
G_X$. 

(ii) $(X,S)$ is indecomposable if and only if the action
$\rho_c$ of the group $G$ is transitive.

(iii) $G_X=A_X$, $\Gamma_1=\Gamma_2=Ker(\rho)\cap A_X$ in Theorem
\ref{Gamma}. Moreover, $i(X)$ generates $G$ and action $\rho_c$ is 
faithful. 

(iv) Starting with any solution $(X, S)$, construct $(X, S')$, 
$S'(x,y)=(\phi(y,x), x)$ for $\phi(y,x)=f_x(g_{f_y^{-1}(x)}(y))$. So
constructed pair $(X, S')$ is a derived
solution, called 
the solution derived from $(X, S)$. Its structure group is
$A_X$. This solution is invariant under the
action of $G_X$, namely $f_z\phi(y,x)=\phi(f_zy,f_zx)$.
\end{theorem}

\subsection{Affine Solutions}

In this sections we recall the characterization of affine solutions that
were studied in great detail in \cite{ESS}, \cite{S}.

Let $X$ be an abelian group.
\begin{definition}
A solution $(X, S)$ is called affine if $S$ is of the
form $S(x,y)=(ax+by+h,cx+dy+t)$, where $a,b,c,d\in End(X)$, $z,t\in X$.
\end{definition}

\begin{lemma}\cite{S}
\label{aff}
(i) Affine solutions are in 1-1 correspondence with 6-tuples
$(q_1,q_2,z,1+s,k,h)\in Aut(X)^4\times X^2$ such that
$$z^2-z(q_1+q_2)+q_1q_2=0,\ q_1q_2=q_2q_1,$$
$$sq_1=q_1s=(1+s)^{-1}sq_2=(1+s)^{-1}q_2s=zs=sz,
$$
$$
sq_1h=(1-q_1)k,\ 
q_1k=zk=(1+s)^{-1}q_2k.
$$

(ii) The correspondence in (i) is given by the formulas
$b=q_1^{-1},\ a=1-zq_1^{-1}, d=1+s-q_1z^{-1}q_2q_1^{-1},
c=b^{-1}((1-d+ad)(1-a)+s)$,
$t=-c(1-a)^{-1}h+k$.

(iii) The affine solution $(X, S)$ is injective 
if and only if $k=0$ and $s=0$ in (i). It is involutive if
and only if, in addition,
$q_1=q_2$. Therefore, injective
affine solutions are in 1-1 correspondence with quadruples 
$(q_1,q_2,z,h)\in Aut(X)^3\times X$ such that $q_1q_2=q_2q_1,\
z^2-z(q_1+q_2)+q_1q_2=0$.
\end{lemma}

We use the above lemma to classify indecomposable affine solutions of
prime order.

{\bf\noindent Example.}

Let $X$ be any set, $f,\ g$ - some permutations from $Aut(X)$. If
$fg=gf$
then $(X, S)$ with $S(x,y)=(g(y),f(x))$ is a solution. 
\begin{definition}
We call the above a
permutation solution.   
\end{definition}

\begin{theorem}
\label{aff_pr}
Let $p$ be a prime number.

1. For each triple $q_1,q_2,h$, $q_1,\ q_2\in \Z_p^*=\Z_{p-1}$, $h\in
\Z_p$,
$q_1\neq q_2$ 
the following are indecomposable nondegenerate affine solutions
of prime order $p$: 

(i) $(x,y)\to (q_1^{-1}y+(1-q_2q_1^{-1})x+h,q_2x-q_1h)$;

(ii) $(x,y)\to (q_1^{-1}y+h,q_2x+(1-q_2q_1^{-1})y-q_2h)$;

(iii) $(x,y)\to (x+h_1,y+h_2)$, $(h_1,h_2)\ne (0,0)$.

2. Any indecomposable nondegenerate affine solution
of prime order $p$ is isomorphic to one
of the above.
\end{theorem}
{\bf\noindent Proof:}

Every affine solution with $|X|=p$ is either injective
or a permutation solution. Indeed, if $s$ of Lemma \ref{aff} is
not
zero  then $q_1=(1+s)^{-1}q_2=z$ and $a=d=0$, so we get a permutation
solution. On the other hand, if $s=0$
then we get that either $q_1=q_2=z=1$ and  we get a permutation
solution  or $k=0$ and we obtain an injective solution.

It is easy to see that indecomposable permutation solutions 
of prime order are affine solutions of type (iii). Indeed, if 
such a solution is given by $(x,y)\to (by,cx)$ then $bc=cb$
and hence $b$ acts transitively on the set of $c$-orbits. 
Therefore, $b$ is cyclic or trivial, which easily implies the 
statement. 

Injective affine solutions are given by triples $(q_1,q_2,z)$ such that
$z^2-z(q_1+q_2)+q_1q_2=0$. This implies that either $z=q_1$ or $z=q_2$ and
we get two solutions introduced in the theorem. Clearly, in the case
$q_1\neq q_2$ these solutions are indecomposable. If $q_1=q_2$ we get that
our solution is involutive. It was proven in \cite{ESS} that involutive
indecomposable solutions are isomorphic to $(\Z_p, S_0)$,
$S_0(x,y)=(y-1,x+1)$, i.e. are permutation solutions. 
The theorem is proved. $\blacksquare$

\section{Indecomposable solutions of prime order}
\begin{theorem}(Main Theorem)
\label{main}
Every indecomposable solution $(X, S)$ of prime order 
$p=|X|$ is affine.
In particular, all
indecomposable solutions of prime order are of the kind considered in
Theorem \ref{aff_pr}. 
\end{theorem}

The rest of the section is occupied by the proof of this theorem. 

In order to prove the main theorem we first classify derived 
solutions with $|X|=p$. The key role in the proof is played 
by the following Lemma. 

\begin{lemma}\label{mainlem} 
An indecomposable derived solution with prime number of elements is 
isomorphic to either

a) $(\Z_p,S)$, where $S(x,y)=(y+1,x)$ or

b) $(\Z_p,S)$, where $S(x,y)=(Ky+(1-K)x,x)$, where $K\in \Z_p,\ K\neq 0,
1$.
\end{lemma}
{\bf\noindent Proof.}

Let $(X, S)$ be derived, $|X|=p$. Then, by Theorem \ref{derived}, there is
a finite group $G$ and $G$-equivariant map $i:X\to G$. Therefore $i(X)$
consists of either one element or $p$ elements. If $|i(X)|=1$ then we get
a permutation solution given by a cyclic permutation ($(X, S)$ has to be
indecomposable), i.e. the solution has the form a). On the other hand, if
$|i(X)|=p$, the map $i$ is injective, and $i(X)$ is a generating conjugacy
class in the group $G$. Also, the group $G=G_X/Z(G_X)$ has no center
(since it must act faithfully on $i(X)$). Thus, by 
Theorem \ref{conjcl} of the Appendix,
$G$ is a subgroup in
an affine group of order $p$. Therefore, our solution has the form b).  
$\blacksquare$

Let $(X, S)$ be an indecomposable solution, $|X|=p$. The $G_X\ltimes A_X$ -
equivariant map $\psi_A:X\to A_X$ gives rise to an equivariant map
$\overline{\psi_A}:X\to A_X/\Gamma_2$, which is either injective or
contracts $X$ into one point. If $X$ is contracted into one point, then
our solution is a permutation solution. It easy to see that an indecomposable 
permutation solution of prime order is affine (see the proof of
Theorem \ref{aff_pr}). 
So in this case the Theorem is proved. 

Suppose now that the map
$\overline{\psi_A}$ is injective. Then the solution $(X, S)$ is injective. 
Since $G_X\ltimes A_X$ acts on $X$, the
action of $G_X$ permutes $A_X$-orbits on $X$. In this way, either the
action of $A_X$ is trivial or transitive. If the action of $A_X$ is
trivial then $A_X$ is abelian, and our solution is involutive
by Theorem \ref{invo}. It was
proven in \cite{ESS} that the only indecomposable involutive solution with
$p$ elements ($p$ is prime) has the form $(\Z_p,S_0)$, where 
$S_0(x,y)=(y+1,x-1)$. 

Assume that the action of $A_X$ on $X$ is transitive and $A_X$ is
not abelian. Then, the solution
$(X, S')$ derived of $(X, S)$ is indecomposable and therefore has the form
$S'(x,y)=(Ky+(1-K)x,x)$, where $X=\Z_p$. 

This implies that 
$S(x,y)=(g_x(y),c(y)x+d(y))$ for some functions $c:\Z_p\to \Z_p$,
$d:\Z_p\to \Z_p$. Indeed, $\phi(y,x)=Ky+(1-K)x$ is $G_X$-invariant, i.e.
$f_z(Ky+(1-K)x)=Kf_z(y)+(1-K)f_z(x)$, thus $f_y^{-1}(x)=c(y)x+d(y)$ by
the following lemma.
\begin{lemma}
If $f:\Z_p\to \Z_p$ satisfies the relation
$f(Ky+(1-K)x)=Kf(y)+(1-K)f(x)$ for some $K\neq 0,1$ and any $x,\ y\in
\Z_p$ then $f$ is affine,
i.e.  there are $c,d\in \Z_p$ such that $f(x)=cx+d$.
\end{lemma}
{\bf\noindent Proof of Lemma:}
Let $\F\subset \Z_p$ be the set of elements $\alpha\in \Z_p$ such that
$f(\alpha x +(1-\alpha)y)=\alpha f(x)+(1-\alpha) f(y)$ for all $x,y\in
\Z_p$. Clearly, $0,1,K,1-K\in \F$. Moreover, if $\alpha,\ \beta,\ \gamma
\in \F$, then $\alpha\beta+(1-\alpha)\gamma\in \F$ since for $x,y$
$$
(\alpha\beta+(1-\alpha)\gamma)x+(1-\alpha\beta-(1-\alpha)\gamma)y=
\alpha(\beta x+(1-\beta)y)+(1-\alpha)(\gamma x +(1-\gamma)y).
$$
Therefore by taking $\gamma=0$ we get that $\alpha\beta\in \F$. In this
way, both $K^{-1}=K^{p-1}$ and $(1-K)^{-1}$ are in $\F$. So, for each
$\alpha,\ \beta\in \F$ the element 
$\alpha+\beta=K \alpha K^{-1}+(1-K)\beta (1-K)^{-1}$ is in $\F$. Since 
$1\in \F$, $\F=\Z_p$ and for any $\alpha\in \Z_p$ 
$f(\alpha x +(1-\alpha) y)=\alpha f(x) + (1-\alpha )f(y)$. If we take
$x=1,\ y=0$, $f(\alpha )=\alpha f(1) + (1-\alpha )f(0)$, i.e. $f$ is
affine. Lemma is proved. $\blacksquare$

Since $f_y^{-1}(x)=c(y)x+d(y)$ is the action of  $G=G_X/\Gamma_1$ on $X$
we can view
$c(y),\ d(y)$ as the functions defined on $G$ that satisfy 
$c(y_1y_2)=c(y_1)c(y_2)$, $d(y_1y_2)=c(y_1)d(y_2) + d(y_1)$.
In particular, $c:G\to \Z_{p-1}$ is a homomorphism of groups.
We would like to show that $c(y)$ is independent of $y$ and $d(y)$ is
affine in $y$. For that we study the properties of group $G$.

Let $H\subset G$ be a subset of $G$ given by $H=\{x^{-1}y|x,y\in
\overline{\psi_G}(X)\}$, where $\overline{\psi_G}:X\to  G$ is coming
from $\psi_G:X\to  G_X$.
\begin{lemma}

(i) $H$ is a subgroup of  group $G$.

(ii) $|H|=p$.
\end{lemma}

{\noindent\bf Proof.}

Let us show that $x^{-1}y$ depends only on $x-y$. For any positive integer
$n$ $x^{-1}y=(x^{-1}y)^nx^{-1}y(x^{-1}y)^{-n}$. Besides, since 
$ztz^{-1}=\phi(t,z)=Kt+(1-K)z$, one has
$$
(x^{-1}y)^nx^{-1}y(x^{-1}y)^{-n}=
(x+nK^{-1}(1-K)(y-x))^{-1}(y+nK^{-1}(1-K)(y-x)).
$$
In this way, $x^{-1}y$
depends only on $x-y$.

 Thus, we have
$x^{-1}yz^{-1}t=x^{-1}yy^{-1}(t+y-z)=x^{-1}(t+y-z)$ implying that $H$ is a
group and the surjective map $j:\Z_p\to H$ given by $j(x-y)=x^{-1}y$ is a
group homomorphism. Since $H$ is not trivial (as $x^{-1}y\neq 1\in G$), 
we have that $j$ is an isomorphism, and thus $|H|=p$.$\blacksquare$

The lemma implies that $c(y)$ is independent of $y$. Indeed, 
$c:G\to \Z_{p-1}$ is a homomorphism that is forced to be trivial
restricted to $H=\Z_p$, i.e. $c(x^{-1}y)=1$, i.e. $c(x)=c(y)$. The map
$d:G\to \Z_p$ restricted to $H$ is a homomorphism from $H$ to
$\Z_p$, therefore $dj:\Z_p\to \Z_p$ is a homomorphism too. This implies
that $d(x^{-1}y)=m(x-y)$ for some integer $m$. If we fix $x$ we get that
$m(x-y)=c(x^{-1})d(y)+d(x^{-1})$, i.e. $d(y)$ is affine in $y$. 

Now when we
know that $f_y^{-1}(x)=cy+d_1x+d_2$ and $\phi(y,x)=Ky+(1-K)x$ we can use
the definition of $\phi(y,x)$ in Theorem \ref{derived}(iv) to conclude
that $g_x(y)=ax+by+h$. Theorem is proved. $\blacksquare$

\section{Appendix: finite groups with trivial center generated 
by a conjugacy class of prime order}

The goal of this section is to prove the following theorem, 
which is used in the proof of Lemma \ref{mainlem}.

\begin{theorem}\label{conjcl} Let $G$ be a finite group 
and $C$ a conjugacy class in $G$ of prime order $p$.  Assume
that the center $Z(G)$ of $G$ is trivial, and 
that $G$ is generated by $C$.  Then
$G$ is a subgroup of the affine group of degree $p$.
\end{theorem} 

{\bf Remark.} In the course of the proof we establish 
some results in the more general case, when $|C|$ is a 
power of a prime. These results may be used in classifying 
indecomposable derived solutions of prime power order. 

We first prove an extension of Burnside's Theorem
(\cite{Go}, 4.3.2) which asserts that in a nonabelian
finite simple group, there are no conjugacy classes
of prime power order.

\begin{lemma}\label{L1} Let $G$ be a finite group and $C$ a conjugacy
class of $G$.  If $\phi:G \rightarrow H$ is a surjective
homomorphism, then $|\phi(C)|$ divides $|C|$.
\end{lemma}

{\bf \noindent Proof:} Let $x \in C$.  Then $|\phi(C)|=|G:H|$ where 
$H = \{ g \in G| x^g \in xK\}$ where $K={\rm ker}\phi$.
Since $H \ge C_G(x)$ and $|C|=|G:C_G(x)|$, the result follows.
$\blacksquare$

The next result extends Burnside's theorem.  
Our proof uses the classification
of outer automorphisms of simple groups. 

\begin{lemma} Let $L$ be a finite nonabelian simple group.  If
$x$ is a nontrivial automorphism of $L$, then
$|L:C_L(x)|$ is not a prime power.
\end{lemma}

{\bf Remark.}
If $x$ is an inner automorphism of $L$, this
is precisely Burnside's theorem. So we may
assume that $x$ is an outer automorphism of $L$.
However, the proof we give will not use this.
 
{\bf \noindent Proof:} 
Let $p^a=|L:C_L(x)|$. 

The list of all subgroups $H$ of a simple group of 
index $p^a$ is given in [Gu].  Aside from a short list
(namely $L(2,11)$, $U(4,2)$, $M_{11}$ and $M_{23}$), it follows
that either $L=A_{p^a}$ or $L=L(d,q)$ with $p^a = (q^d-1)/(q-1)$.

In the first four cases, one checks directly that no automorphism
of $L$ centralizes the appropriate subgroup.  If $L=A_{p^a}$,
then $H=A_{p^a-1}$ and it has trivial centralizer in
$S_{p^a}= {\rm Aut(L)}$.  

Finally, consider the case $L=L(d,q)$.  Then $H$ is either
the stabilizer of a $1$-space or hyperplane in the natural
$d$-dimensional module $V$ for $L$.   Consider
the full automorphism group $J$ of
$L$ which is  generated by the group of semilinear
automorphisms of $V$ ($P\Gamma L(d,q)$) and the transpose
inverse map (for $d > 2$).  This latter does not fix
the conjugacy class of $H$ and so does not normalize $H$.
Thus, the normalizer of $H$ in $J$ is the subgroup of semilinear
transformations fixing the $1$-space (or hyperplane).
It is elementary to see that this subgroup has no center
and so is the not the centralizer of any automorphism.
$\blacksquare$ 

\begin{lemma}  Let $G$ be a finite group and $x \in G$.
Let $C=x^G$ be the conjugacy class of $x$.
Let $N$ be the normal subgroup of $G$ generated
by $C$.  If $C$ has prime power order, then $N$
is solvable.
\end{lemma}

{\bf \noindent Proof:} 
 Assume that $G$ is a minimal counterexample
to the theorem.  Let $A$ be a minimal normal subgroup of 
$G$ contained in $N$.  By Lemma 5,
the image of $C$ in $G/A$ also has prime power order
and so by minimality, $N/A$ is solvable.  So if
$A$ is solvable, the result follows.  

So we may assume that $A$ is a direct product of isomorphic
nonabelian simple groups.  If $x$ commutes with $A$, then
so does $C$ and so $N$.  This implies that $N$ is abelian,
a contradiction.  

Let $L$ be a direct factor of $N$. Let $|C|=p^a$
with $p$ prime. It follows that $x$ commutes with
a Sylow $r$-subgroup of $G$ for every prime $r \ne p$.
In particular, choose $r \ne p$ dividing the order
of $A$.  Since any Sylow $r$-subgroup of $A$ intersects
each direct factor of $A$ and since $x$ permutes the simple direct
factors of $A$, it follows that $x$ normalizes each
direct factor of $A$.  Let $L$ denote a simple direct
factor of $A$.  Since $x$ does not commute with $A$,
we may choose $L$ so that $x$ does not centralize $L$.
By the previous lemma, $|L:C_L(x)|$ is not a prime power
and so neither is $|A:C_A(x)|$.  Since $A$ is normal
in $G$, $|A:C_A(x)|$ divides $|C|$, a contradiction.
$\blacksquare$

With the previous result at hand, we can 
pin down the structure of the normal subgroup
generated by a conjugacy class of prime power order.

Let $O_p(H)$ denote the maximal normal $p$-subgroup of 
a finite group $H$.

\begin{theorem} Let $G$ be a finite group and $C$ a conjugacy class
of $G$ of order $p^a$ with $p$ prime.  Let $N=\langle C \rangle$.
Then  $N/O_p(N)$ is abelian.  In particular, if $G=N$, then
$G/O_p(G)$ is cyclic.
\end{theorem}

{\bf \noindent Proof:} 
  The last statement follows from the first one.

Let $A$ be a minimal normal subgroup of $G$ contained in $N$.
Since $N$ is solvable, it follows that $A$ is an elementary
abelian $r$-group for some prime $r$.

If $r=p$, then the result follows by induction (considering $G/A$
and the image of $C$ in $G/A$).  If $r \ne p$, then $A$ is contained
in every Sylow $r$-subgroup of $G$.  Since $x \in C$ implies that
$C_G(x)$ contains some Sylow $r$-subgroup, it follows that
$C$ commutes with $A$.  Thus, $A \subset Z(N)$.   By considering
$G/A$, we see that either $O_p(N/A) \ne 1$ or $N/A$ is abelian. 
 
Suppose that $O_p(N/A) \ne 1$.  Let $B \subset N$ with $B/A=O_p(N/A)$.
Then $B/Z(B)$ is a $p$-group and so $B$ is nilpotent.  Thus,
$O_p(B) \ne 1$ and since $B$ is normal in $G$, it follows
that $O_p(N) \ne 1$, a case already dealt with.

Suppose that $N/A$ is abelian.  Then $N/Z(N)$ is abelian and so
$N$ is nilpotent.  Since we may assume that $O_p(N)=1$, $N$ is
a $p'$-group.  Let $x \in C$.  Then $C_G(x)$ contains a Sylow $q$-subgroup
for every prime $q \ne p$ and so contains the normal (in $G$) Sylow
$q$-subgroup of $N$.  Thus, $x \in Z(N)$ and so $N$ is abelian.
$\blacksquare$

{\bf Proof of the theorem.}
  Map $G$ into $S_p$ by letting $G$ act on
$C$ by conjugation.  The kernel of this map
is $C_G(C)=C_G(G)=Z(G)=1$.  So this is an embedding.

By the previous result, $G/O_p(G)$ is cyclic.  Since
$O_p(G) \ne 1$ (or $G$ would not be transitive), it follows
that $G$ is contained in the normalizer of a cyclic subgroup
of order $p$ as desired.

{\bf Acknowledgments.}
The work of P.E. was partially supported 
by the NSF grant DMS-9700477, and was partly done 
when he was a CMI prize fellow. P.E. thanks IHES for hospitality. 
The work of R.G. was partially supported by the NSF grant 
DMS-9970305.

\end{document}